\documentclass[english]{article}
\usepackage[T1]{fontenc}
\usepackage[latin9]{inputenc}
\usepackage{amsmath}
\usepackage{graphicx}

\usepackage{babel}

\begin{document}

\title{Algorithmic proof of Barnette's
Conjecture}

\author{I. Cahit \\Near East University\\\texttt{email: icahit@gmail.com}}
\date{}
\maketitle

\begin{abstract}
In this paper we have given an algorithmic proof of an long standing Barnette's conjecture (1969) that every $3$-connected bipartite cubic planar graph is hamiltonian. Our method is quite different than the known approaches and it rely on the operation of opening disjoint \emph{chambers}, by using spiral-chain like movement of the outer-cycle elastic-sticky edges of the cubic planar graph. In fact we have shown that in hamiltonicity of Barnette-graph a single-chamber or double-chamber with a bridge face is enough to transform the problem into finding specific Hamilton path in the cubic bipartite planar graph reduced. In the last part of the paper we have demonstrated that, if the given cubic planar graph is non-hamiltonian, then the algorithm which constructs spiral-chain (or double-spiral chain) like chamber shows that except one vertex there exists $(n-1)$-vertex cycle.

\end{abstract}

\vspace{1cm}

\textbf{\large {1 Introduction}}
\vspace{0.3cm}

Spanning cycle of dodecahedron is the origin of the famous Hamiltonian cycle problem in graphs. Next is the Tait's "conjecture" of hamiltonicity of cubic planar graphs which has been shown to be wrong by Tutte is another wave of stimulation of research area [1],[7]. The best characterization of Hamiltonian graphs was given in 1972 by Bondy and Chvátal theorem which generalizes earlier results by Dirac and Ore [2].

\vspace{0.2cm}

\emph{\textbf{Theorem 1} (Bondy and Chvátal). A graph is Hamiltonian iff its closure is Hamiltonian. }

\vspace{0.2cm}

Given a graph $G$ with $n$ vertices the closure $cl(G)$ is uniquely constructed from $G$ by successively adding for all nonadjacent pairs of vertices $u$ and $v$ with $deg(u)+deg(v)\geq n$ the new edge $uv$.

In general hamiltonian cycle problem in graphs is NP-complete, and remain NP-complete for perfect graphs, planar bipartite graphs, grid graphs, $3$-connected planar graphs [2]. However polynomial algorithm has been given by Gihiba and Nishizeki (1989) for $4$-connected planar graphs [3],[4],[5],[6]. Hence our algorithm is important since it shows that hamiltonicity of Barnette graph in linear time.

Barnette has made the following conjecture in 1969 [8]:

\vspace{0.2cm}

\emph{\textbf{Conjecture 1,(Barnette),1969)}. Every graph that is $3$-connected, $3$-regular, bipartite and planar has a hamiltonian cycle.
}
\vspace{0.2cm}

\begin{figure}
\centering
\includegraphics[scale=0.3]{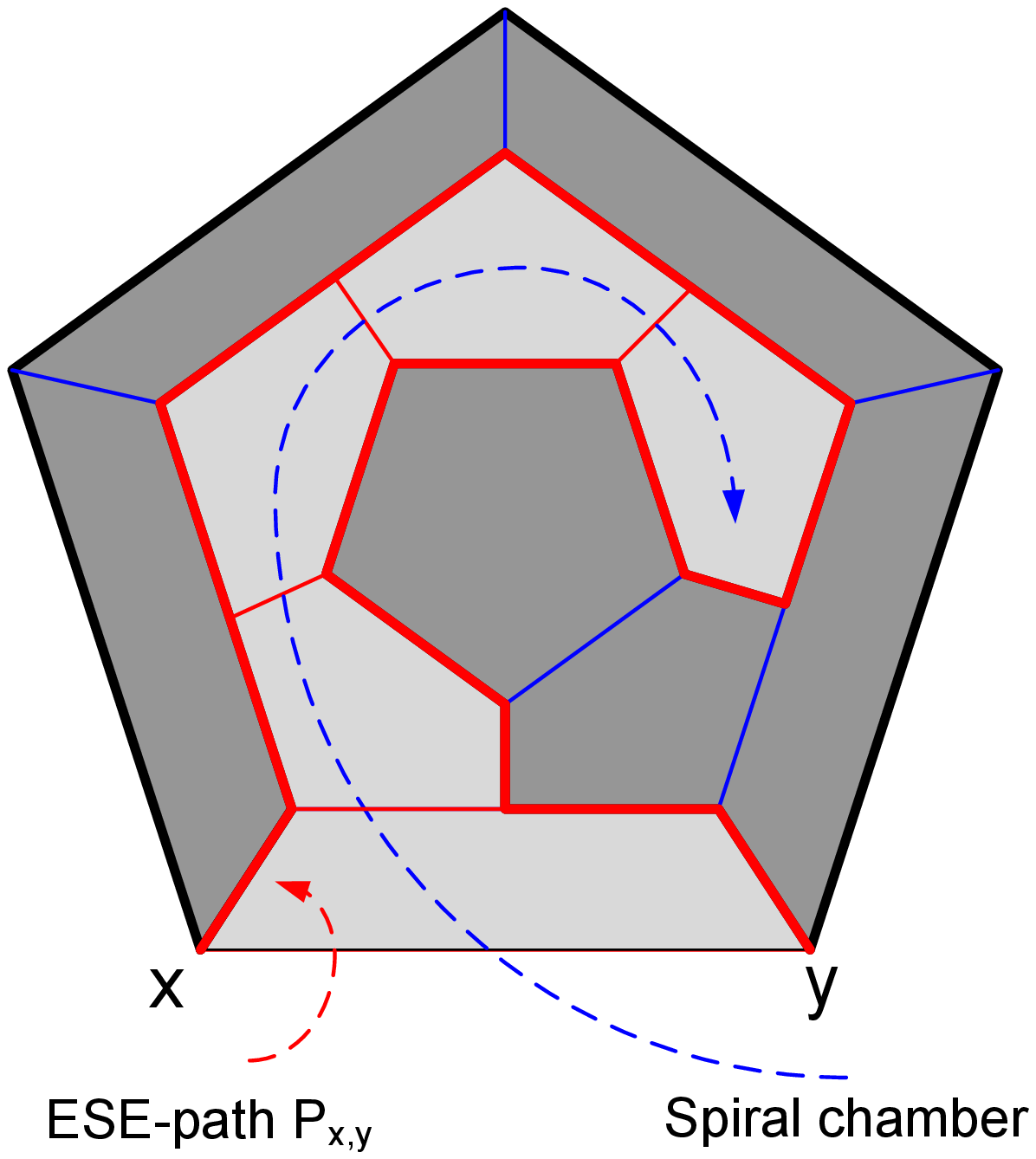}
\caption{Hamilton cycle of dodecahedron, spiral chamber and ESE path.}
\end{figure}

Any graph satisfying the conditions of Conjecture 1 is called Barnette-graph. An excellent survey together some new ideas on Barnette graphs has been given by Luis de la Torre [4]. In fact algorithmic proof given in this paper, is related with an stronger conjecture than Barnette's conjecture which is based on hamiltonian cycles of a list of  Tutte embeddings of Barnette graphs from $8$ to $16$ vertices (see Appendix A [4]). We will give also an argument to rule out the possibility of existence of Tutte's fragments in the Barnette graphs. Similar results have been obtained using a different approach by Kim and Lee in [9] . First Temperley-Lieb algebras have been generalized to \textrm{sl(3,C)} web spaces. Since a cubic bipartite planar graph with suitable directions on edges is a web, the quantum \textrm{sl(3)} invariants naturally extend to all cubic bipartite graph. They completely classify cubic bipartite planar graphs as a connected sum of primes webs and provide a method to find all prime webs and exhibit all prime web up to $20$ vertices. Goodey showed the conjecture holds when all faces of the graph have either $4$ or $6$ sides [10],[11]. Feder and Subi generalize this by showing that when the faces of such graph are $3$-colored, with adjacent faces having different colors, if two of the three color classes contain only faces with either $4$ or $6$, then the conjecture holds [12]. Kelmans has shown the following important theorem which is equivalent to Conjecture 1 [13]:

\vspace{0.2cm}

\emph{\textbf{Theorem 2 (Kelmans)}. (a) For every bipartite, cubic, $3$-connected and planar graph $G$ and for every edges $a,b$ of $G$, belonging to the same facial face of $G$, there is a hamiltonian cycle in $G$ containing $a$ and avoiding $b$.\\
(b)For every bipartite, cubic, $3$-connected and planar graph $G$ and for every edges $a,b$ of $G$, belonging to the same facial face of $G$, there is a hamiltonian cycle in $G$ containing both $a$ and $b$.}

\vspace{0.2cm}

\begin{figure}
\centering
\includegraphics[scale=0.5]{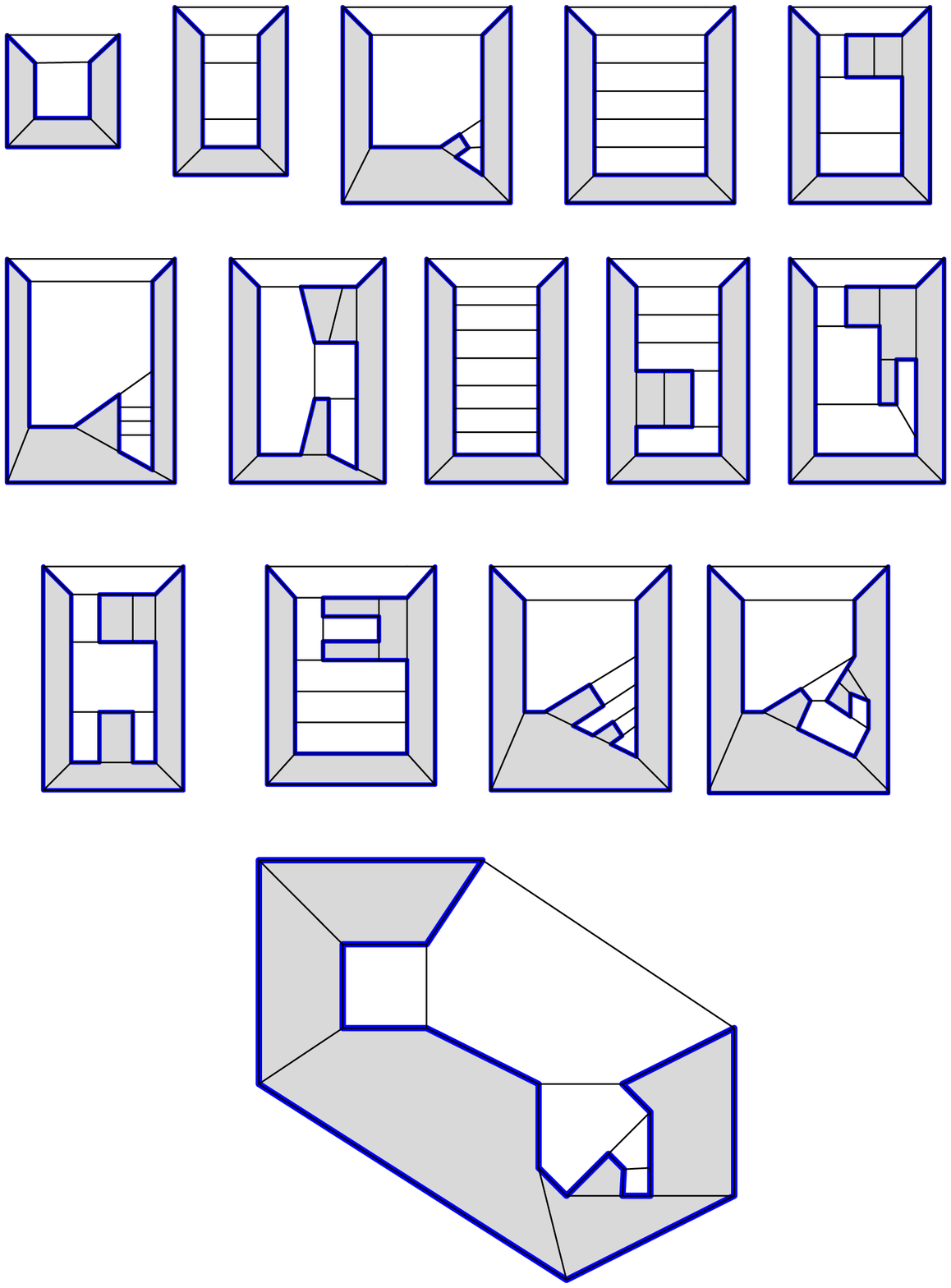}
\caption{Hamilton cycles in Barnette graphs.}
\end{figure}

Hertel has given stronger than Theorem 2 [17].\\

\vspace{0.2cm}

\emph{\textbf{Theorem 3 (Hertel)}. Barnette's conjecture holds if and only if any arbitrary path $P$ of length $3$ that lies on a face in a Barnette graph is a hamiltonian cycle which passes through the middle edge in $P$ and avoids its leading and trailing edges.}

\vspace{0.2cm}

Holton et.al., have shown that $3$ connected cubic graphs with fewer than $66$ vertices are hamiltonian [14] together the relations of $3$-cut and essential $4$-cut with the possible smallest non-hamiltonian graphs.  Aldler \emph{et.al.,} have announced that through a computer search Conjecture 1 is true at most for $84$ vertices [15].

\vspace{0.5cm}

\textbf{\large {2  Algorithmic proof of Conjecture 1}}
\vspace{0.3cm}

Let $G$ be denote a cubic, bipartite planar graph with $n$ vertices. Assume that $G$ drawn suitably in the plane that no edges crosses each other. $C_{o}$ denotes outer-cycle of $G$, where $|C_{o}|\geq 4$.  By $H$ we denote a hamiltonian cycle which passes through all vertices of $G$ such that its edge set partitioned into two subsets
$E(H)=H_{o}\cup H_{i}$, where $H_{o}=\{h_{o,j}\in E(C_{o}),j=1,2,...,k-1\}$ and $H_{i}=\{h_{i,j}\notin E(C_{o}),j=1,2,...,m\}$, $n=k+m-1$. Hence the edge set of $C_{o}=\{h_{o,1},h_{o,2},...,h_{o,k-1}\}\cup \{d_{e}\}$ where subscript $e$ indicates the \emph{entrance} edge of the outer-cycle which is not in $H$. Hence the set of edges of $G$ can be expressed as\\

$E(G)=H_{o}\cup H_{i}\cup D_{o}\cup D_{i}\cup \{d_{e} \}$\\

where the set $D_{o}$ denotes the door-edges remain outside of the region bounded by the hamiltonian cycle $H$ and the set $D_{i}$ denotes the door-edges remain inside the region bounded by the hamiltonian cycle $H$ and  $\{d_{e}\}$ denotes the entrance door-edge. We also note that the number of entrance door-edges may be more than one for an single-chamber. For example double-spiral shape hamiltonian cycle $H$ shown in Figure 5 ($10_{4}$) has two entrance doors $d_{e1}$ and $d_{e2}$.

\vspace{0.2cm}

\emph{\textbf{Definition 1}. The cycle $C_{c}= \{H_{i}\}\cup \{d_{e}\}_{j}$ is called the chamber-cycle induced by the hamiltonian cycle $H$ of $G$.}

\vspace{0.2cm}

If for an hamiltonian cycle $H$ of $G$ there is only one chamber-cycle $C_{c}$ as above we say single-chambered $H$ (see Figures 1) otherwise we call it multi-chambered $H$. It is easy to see that for any hamiltonian cycle $H$ of $G$ no two door-edges $d_{i}$ and $d_{j}$ are adjacent.

In Figure 2 we have shown single-chamber hamiltonian cycles of all Barnette graphs from $8$ to $16$ vertices. In Figure 4 we also give single-chamber hamiltonian cycles of all prime webs up to $20$ vertices [9]. This gives us encourage to state and prove the following:

\vspace{0.2cm}

\emph{\textbf{Conjecture 2}. All Barnette graphs with at most one $3$-cut have single-chamber hamiltonian cycles.}

\vspace{0.2cm}

Clearly the restriction of single-chamber hamiltonian cycle $H$ in $G$ makes the Conjecture 2 easier to prove or disprove than the Conjecture 1. That is, right from the beginning we assume that all outer-edges (except $d_{e}$) of $H_{o}$ are readily in the hamiltonian cycle $H$. Hence if $x$ and $y$ are the end points of the entrance-edge $d_{e}$ the hamiltonian cycle problem would reduce to find an hamiltonian path $P_{H}(x,y)$ in the subgraph $G_{1}=G\setminus\{H_{o}\}$. In the Algorithm below hamiltonian path is constructed step-by-step by stretching the entrance-edge $d_{e}$ onto the edges of the chamber. We will call this operation as adding elastic-sticky edge.

\vspace{0.3cm}

\textbf{\textbf  {2.1  A possible threat to Conjectures 1 and 2}}
\vspace{0.3cm}

Tutte has given a counterexample to Tait's conjecture that all $3$-connected cubic planar graphs have hamiltonian cycles. The main element of the counter-example now is known as Tutte's fragment shown in Figure 3(a) with three critical vertices $x,y,z$ on the corners of the fragment. A sub-hamiltonian paths $P_{H}(i,j)$ only exists if $i\in \{x,y\}$ and $j=z$. Now if one can construct a fragment with three corners by using only \emph{even} cycles that would be a counter-example both for Conjectures 1 and 2. Closest constructions using only cycles of lengths $4$ and $6$ is shown in Figure 3(b) and (c) with $13$ vertices and fortunately they fail. This is true in general, since for any sub-hamiltonian path around an even cycle no vertex of an even cycle can be left unvisited or end-vertex of the sub-hamiltonian path. This observation is equivalent, in the Algorithm 1, that no two door-edges $d_{i}$ and $d_{j}$ would adjacent in the chamber cycle $C_{c}$. This is always possible since all faces in $G$ are even. This is clearly seen, then algorithm applied for non-hamiltonian planar graphs (see Figure 5).

\begin{figure}
\centering
\includegraphics[scale=0.4]{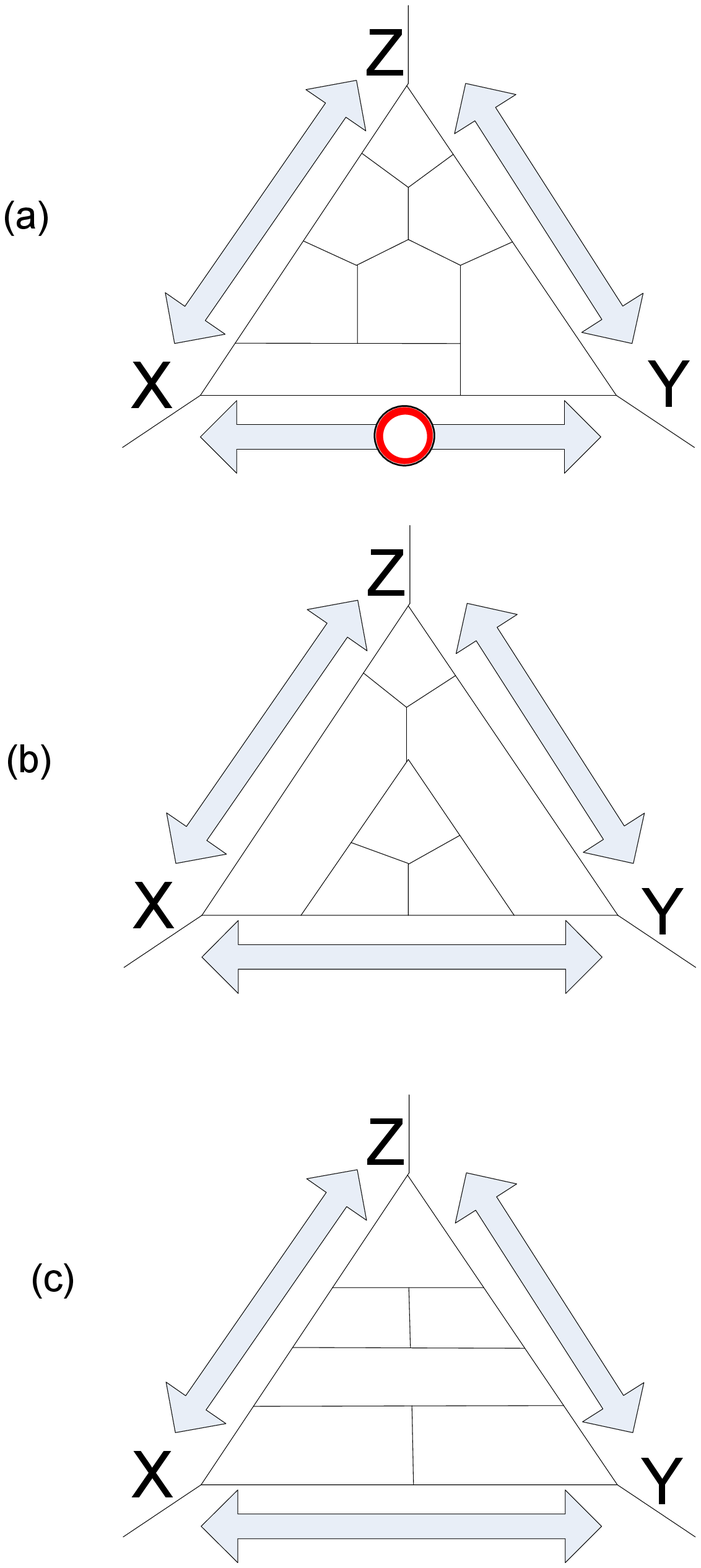}
\caption{(a) The Tutte's fragment, (b),(c) unsuccessful bipartite fragments.}
\end{figure}
\begin{figure}
\centering
\includegraphics[scale=0.5]{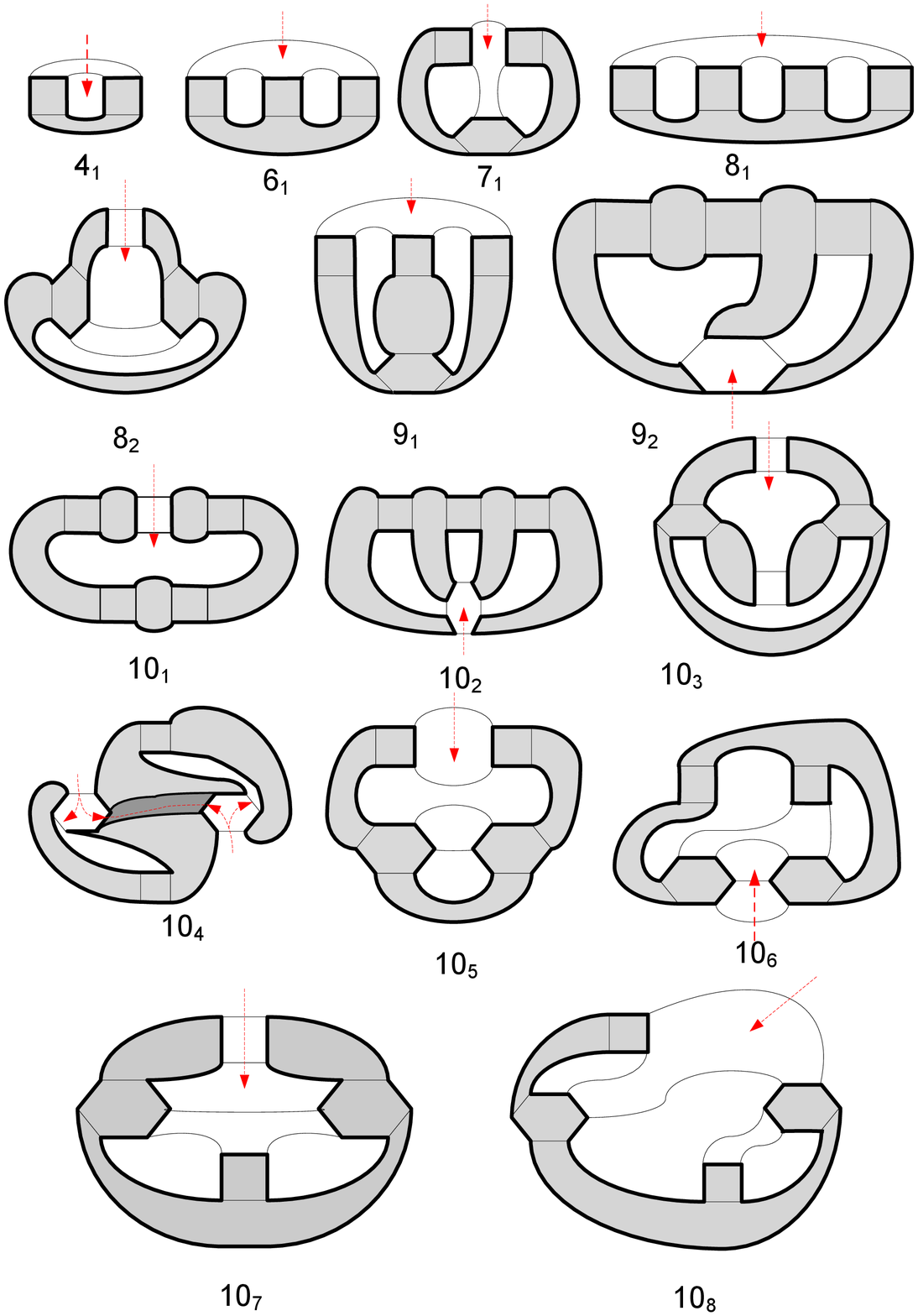}
\caption{Hamiltonian cycles under the quantum sl(3) invariants of Barnette graphs up to $20$ vertices.}
\end{figure}

\begin{figure}
\centering
\includegraphics[scale=0.5]{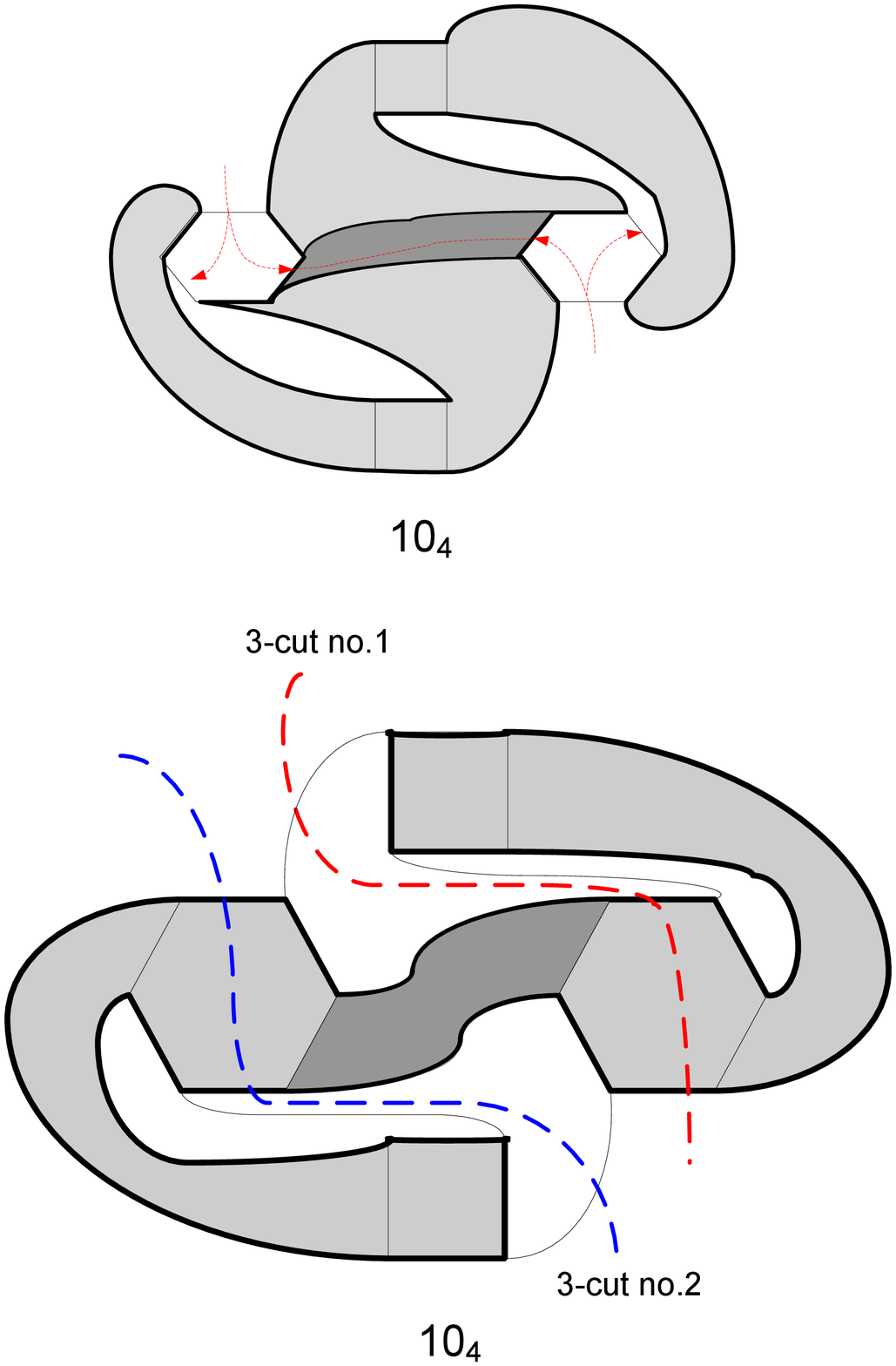}
\caption{Hamiltonian cycles in double-spiral shape of web $10_{4}$ and two edge disjoints $3$-cuts, the bridge face is shown in dark-gray.}
\end{figure}

\vspace{0.3cm}
\textbf{\textbf  {2.2 The algorithm carve-cubic-planar}}
\vspace{0.3cm}

Let us start with a useful Lemma.

\vspace{0.2cm}

\emph{\textbf{Lemma 1}. Let $G$ be a Barnette graph with a $3$-cut $\{a,b,c\}$, $a,b,c\in E(G)$. That is $G=G_{1}\cup G_{2}\cup \{a,b,c\}$. Then in any single-chamber hamiltonian cycle $H$ the entrance-door edge $d_{e}\notin G_{1}$ or $G_{2}$ .}\\

\vspace{0.2cm}

\emph{Proof:} If the edges $a,b,c$ are the $3$-cut, where $a$ and $c$ are outer-cycle edges, then hamiltonian cycle $H$ must contains both $a$ and $b$ or both $b$ and $c$. Either case implies another chamber by the entrance-door edge $d_{e}=c$ or $d_{e}=a$.

\emph{\textbf{Algorithm }(Carve-Cubic-Planar):}\\
\emph{Step 1:} \emph{(Initial Chamber)}.\\
Let $G$ be a $3$-connected, bipartite cubic planar graph. First select a suitable outer-edge (see Lemma 1) for the entrance door-edge $d_{e}$. Hence outer-edges of $G$ is $E_{o}=\{d_{e},h_{o,1},h_{o,2},...,h_{o,k}\}$, where $k+1$ is even. Initially the entrance door-edge defines a facial cycle (face) $C_{c,1}=\{d_{e},e_{in,1},e_{in,2},...,e_{in,r}\}$. Since $|C_{c,1}|$ is even we can rewrite its edges as $C_{c,j}=\{d_{e},h_{in,1},e_{in,2},h_{in,3},...,h_{in,r}\}$. That is $e_{in,j}=h_{in,j},j=1,3,...,r$ becomes subset of internal hamiltonian edges and $e_{in,j}=d_{in,j},j=2,4,...,r-1$ becomes internal door-edges. Hence $H=\{H_{o}\cup H_{in,c1}\}$ where $H_{in,c1}$ is the set of internal hamiltonian cycle edges of the chamber $C_{c}$ defined by $d_{e}$. Similarly let $D_{in,c1}$ be the set of door-edges defined by $d_{e}$.\\
\emph{Step 2:} \emph{(Knock-the door and enter)}.\\
Repeat Step 1 for each door-edge $d_{i}\in D_{in,c1}$. If the face (cycle) defined by door-edge $d_{i}$ contains an edge that share a cycle from the set $H_{o}$ then we put the edge into the set of internal hamiltonian edges. That is $H=\{H_{o}\cup H_{in,c1}\cup H_{in,d_{i}}\}$, $d_{i}\in D_{in,c1}$. If the door-edge $d_{i}$ defines a face (cycle) $C_{i}$ has an edge $e$ which share a cycle $C_{j}$ with an outer-hamiltonian edge and with another door-edge:\\
$C_{i}=\{e,d_{i},...\}$, $C_{j}=\{e,h_{oj},d_{j},...\}$, $e=\{C_{i}\cap C_{j}\},h_{oj}\in H_{o}, d_{i}\neq d_{j}$
Then put the edge $e$ and door-edge into the set $H_{i}$. That is we call the cycle (face) as the bridge-face in the hamiltonian cycle $H$ (see Figure 5). This situation arises when $G$ has two edge-disjoint $3$-cuts . Algorithm continue  from the door edge  $d_{j}$.
If $|H|=n$ then we have entered all faces through the door-edges and a hamiltonian cycle has been found. Otherwise we repeat Step 1 for the other door-edges in the other levels. Note here that we have not selected adjacent door-edges.\\
Illustration of algorithm is shown in Figure 6.

\vspace{0.2cm}

\emph{\textbf{Theorem 4.} Let $G$ be any cubic, $3$-connected, bipartite planar graph $G$. Then the Algorithm "Carve-Cubic-Planar" always terminate with an hamiltonian cycle $H$ of $G$.}

\vspace{0.2cm}

\emph{Proof.} Let us assume that algorithm CCP has not produced a hamiltonian cycle $H$. Then there must be a vertex $v_{x}\notin H$ and $v_{x}$ must be exactly in three facial cycles $C_{1},C_{2}$ and $C_{3}$. Without loss of generality assume that step "knock-the door and enter" has been performed for $C_{1}$ before $C_{2}$ and $C_{3}$. Then there must be two vertices $v_{y}$ and $v_{z}$ of $C_{1}$ such that $(v_{x},v_{y}),(v_{x},v_{z})\in C_{1}$. Then we see that both edges would be door-edges and the cycle $C_{1}$ is odd.

\begin{figure}
\centering
\includegraphics[scale=0.5]{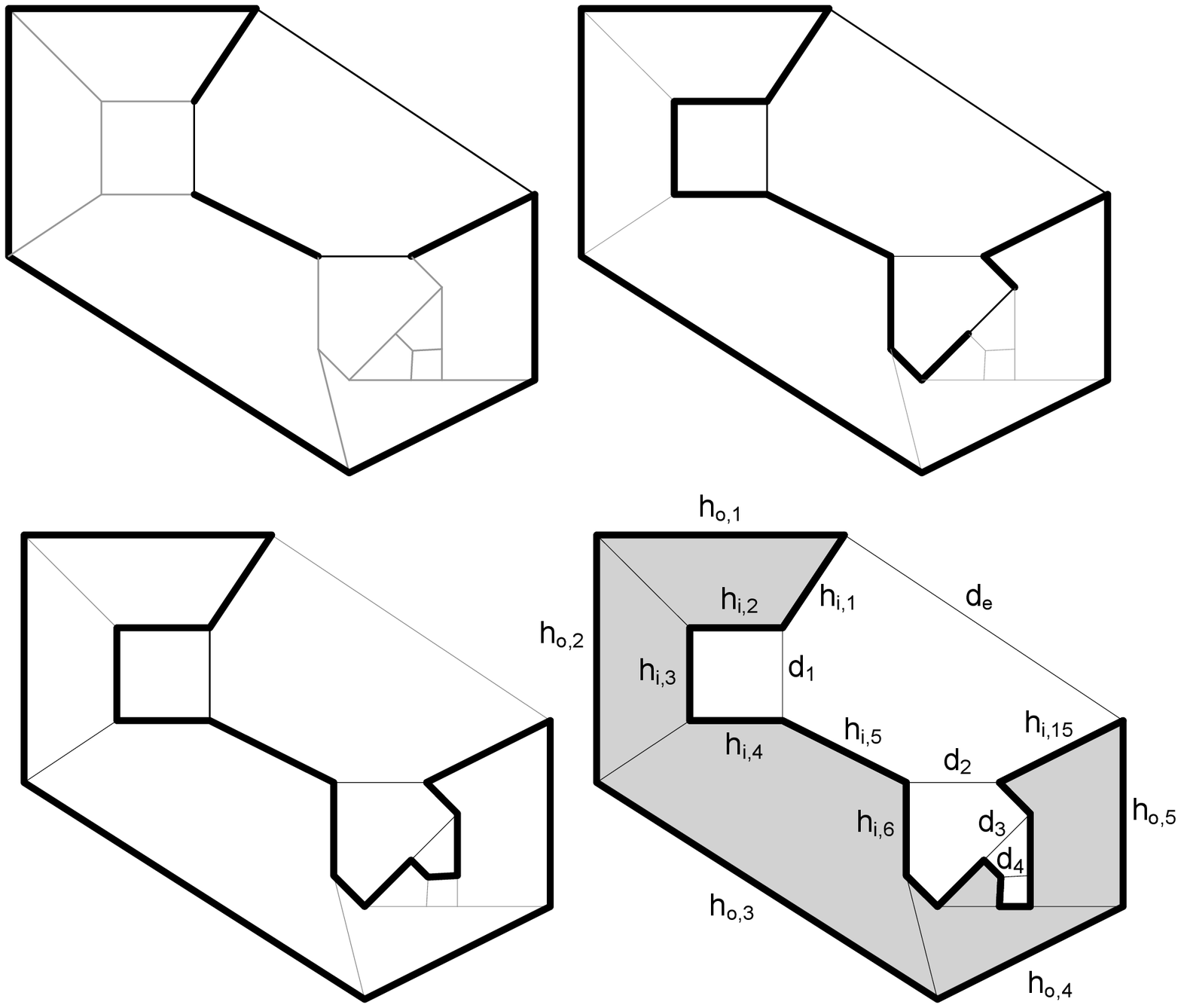}
\caption{Implementation of Algorithm Carve-Cubic-Planar.}
\end{figure}

\vspace{0.3cm}
\textbf{\large {3 Non-hamiltonian $3$-connected cubic planar graphs}}
\vspace{0.3cm}

Holton et.al., have shown that all $3$-connected cubic planar graphs on $36$ or fewer vertices are hamiltonian and the only non-hamiltonian examples on $38$ vertices which are not cyclically $4$-connected are the six graphs which have been found by Lederberg, Barnette and Bos\'{a}k [16]. We have shown non-hamiltonian cubic planar graphs with $42,46$ and $44$ vertices in Figure 7  [16] together with cycles of length $n-1$. As shown in Figures 7(a) and (b), if we choose right-door edges in the chamber cycles in the algorithm the resulting longest cycles are in the shape of spiral $S$. We can alternatively select two entrance door edges  symmetrically, the algorithm again results an $(n-1)$-vertex cycle in the forum of a double-spiral $S_{1}$ and $S_{2}$.

\vspace{0.2cm}

\emph{\textbf{Theorem 5}. For every non-hamiltonian $3$-connected cubic planar graph, Algorithm 1 terminates with a cycle of length $n-1$.}

\vspace{0.2cm}

\begin{figure}
\centering
\includegraphics[scale=0.5]{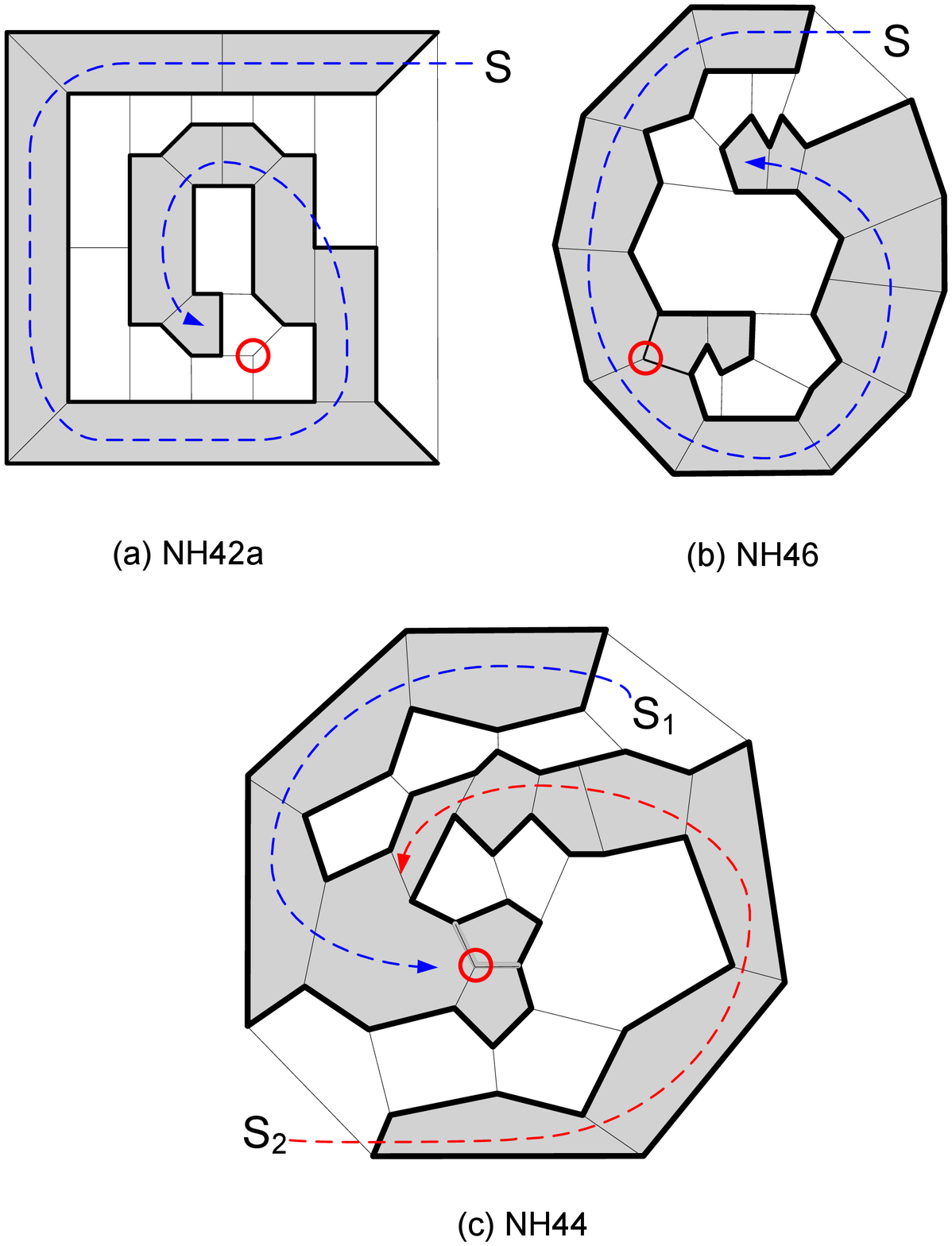}
\caption{Some non-hamiltonian cubic planar graphs.}
\end{figure}

\vspace{0.3cm}
\textbf{\large {4   Concluding remarks}}
\vspace{0.3cm}

In this paper we have given an algorithmic proof of Barnette's conjecture that all $3$-connected bipartite cubic planar graphs is hamiltonian. The algorithm given for this, which we called "carve-cubic-planar" has some interesting features: First of all it delete the edges (door-edges) in the expanding chamber cycle that will not be in the hamiltonian cycle. Secondly by selecting specific entrance door-edge for the chamber, all outer-edges of the graph becomes edges of the seeking hamiltonian cycle and hence the problem reduced of finding hamiltonian path in the remaining cubic planar graph. Lastly hamiltonian cycle as seen a chain of faces by the algorithm would look like a spiral. This has been particularly demonstrated through the examples of non-hamiltonicity of $3$-connected cubic planar graphs. Since spiral-chain coloring algorithm has been unified for the solution of several graph coloring problems [18], it may as well be used in the solution of other problems related hamiltonian cycles of planar graphs in general [19].


\newpage
\begin{thebibliography}{22}

\bibitem[1]{key-1}W. T. Tutte, On hamiltonian circuits, \emph{J. London Math. Soc.},21(1946),98-101.
\bibitem[2]{key-2}J. A. Bondy and U.S.R. Murty, \emph{Graph Theory with Applications}, MacMillan and Co.,London 1976.
\bibitem[3]{key-3}M.R.Garey and D.S.Johnson,\emph{Computers and Intractability: A Guide to the Theory of NP-Completeness},W.H. Freeman and Co., New York, NY, 1990.
\bibitem[4]{key-4}L. de la Torre, \emph{Investigations of Barnette's Graphs}, Undergraduate Thesis, Dept. of Math., Univ. of California, Davis, 2005.
\bibitem[5]{key-5} M.R. Garey, D.S. Johnson and R.E. Tarjan, The planar Hamiltonian circuit problem is NP-complete, \emph{SIAM J. of Computing}, 5 (1976).
\bibitem[6]{key-6} T. Akiyama, T. Nishizeki, N. Saito, NP-completeness of the Hamiltonian cycle problem for bipartite graphs, \emph{J. of Info. Processing} ,(1980).
\bibitem[7]{key-7}R. Greenlaw, R. Petreschi, Cubic graphs, \emph{ACM Computing Surveys}, 27(4), December 1995, 471-495.
\bibitem[8]{key-8}D. Barnette, Conjecture 5, \emph{Recent Progress in Combinatorics} (Ed. W.T. Tutte), Acad. Press, New York, (1969),343.
\bibitem[9]{key-9}D. Kim and J. Lee, The quantum sl(3) invariants of cubic bipartite planar graphs, \emph{J. of Knot Theory and its Ramifications}, 17(3), (2008),361-375.
\bibitem[10]{key-10}P.R. Goodey, Hamiltonian circuits in polytopes with even sides, \emph{Israel J. Math.}, 22 (1975), 52-56.
\bibitem[11]{key-11}P.R. Goodey, A class of Hamiltonian polytopes, \emph{J. Graph Theory}, 1(1977), 181-185.
\bibitem[12]{key-12} T. Feder and C. Subi, On Barnette's conjecture, \emph{ECCC Report}, TR05-005, 2005.
\bibitem[13]{key-13}Kelmans, On Hamiltonian cycles in bipartite cubic $3$-connected planar graphs, Tech. Rep. Rutcor. Res. Rep., 26 June 2003.
\bibitem[14]{key-14}D. A. Holton, B. Manvel, B. D. McKay, Hamiltonian cycles in cubic $3$-connected bipartite planar graphs, \emph{J. Combin. Theory} (Series B),38(3), jUNE 1985,279-287.
\bibitem[15]{key-15}R. Aldler, G. Brinkmann and B. McKay, Research Announcement, December 2002.
\bibitem[16]{key-16}D. A. Holton, B.D. McKay, The smallest non-hamiltonian $3$-connected cubic planar graphs, \emph{J. Combin. Theory} (Series B), 45(3), December 1988, 305-318.
\bibitem[17]{key-17}A. Hertel, A survey and strengthening of Barnette's conjecture, April 2005, preprint.
\bibitem[18]{key-18}I. Cahit, A unified spiral chain coloring algorithm for planar graphs, in \emph{Structural Analysis of Complex Networks}, D. Matthias (Ed.) 2009,Birkhäuser-Boston.
\bibitem[19]{key-19}J.A. Bondy, U.S.R. Murty, \emph{Graph Theory: An Advanced Course}, Springer 2007, 651 pages.





\end{thebibliography}
\end{document}